\documentclass[journal]{IEEEtran}
\usepackage[utf8]{inputenc}
\usepackage{textcomp}
\usepackage{graphicx}
\usepackage{amsmath}
\usepackage[version=4]{mhchem}
\usepackage{siunitx}
\usepackage{longtable,tabularx}
\begin{document}

\title{Optimal Disturbance Attenuation Approach with Measurement Feedback to Missile Guidance}

\author{Barak Or, Joseph Z. Ben-Asher and~Isaac Yaesh
\thanks{B. Or and  Joseph Z. Ben-Asher are with the Technion, Israel Institute of Technology,
IL, 32000 e-mail: (barakorr@gmail.com).}
\thanks{I. Yaesh is with ELBIT SYSTEMS LTD}
}

\markboth{}%
{Shell \MakeLowercase{\textit{et al.}}: Bare Demo of IEEEtran.cls for IEEE Journals}

\maketitle

\begin{abstract}
Pursuit-evasion differential games  using the Disturbance Attenuation approach are revisited. Under this approach, the pursuer actions are considered to be control actions, whereas
all external actions, such as target maneuvers, measurement errors and initial position uncertainties, are considered to be disturbances.
Two open issues have been addressed, namely the effect of noise on the control gains, and the effect of trajectory shaping on the solution. These issues are closely related to the question of the best choice for the disturbance attenuation ratio. Detailed analyses are performed for two simple pursuit-evasion cases: a Simple Boat Guidance Problem and Missile Guidance Engagement.
\end{abstract}

\begin{IEEEkeywords}
Missile Guidance, Disturbance Attenuation, Pursuit-evasion games, Integration of Guidance /  Estimation.
\end{IEEEkeywords}

\section{Introduction}
\IEEEPARstart{T}{he} seminal work of  Ref. \cite{isaacs1999differential} was the first to formulate the pursuit-evasion problem  as a zero-sum differential game between the pursuer and the evader. Many researchers have followed Isaac's footsteps. In particular, Linear Quadratic Differential Games, (LQDG) \cite{bryson1975applied}, \cite{ben1998advances} have attracted attention.  Differential games where measurement noises are considered as the evader's co-players have also been investigated in Ref. \cite{speyer1976adaptive} and \cite{ben2017games}, under the so-called disturbance attenuation (DA) approach. In these publications a simple model of pursuer was realized, and the relation between the noise magnitude and the control gain was obtained. Quite unexpectedly, the control gain grows as the noise grows. The first objective of the present paper is to investigate this phenomenon.
The LQDG formulation enables adding an additional term to the cost function, namely a trajectory shaping term. It was shown in \cite{ben2004trajectory} that this term may improve the interception performance. Additionally, it may lead to some improvement in the estimation performance \cite{oshman1999optimization}. In previous work with the DA approach(e.g. Ref. \cite{speyer1976adaptive}) no trajectory shaping terms were considered. In this case the estimator equation reduces to a simple Kalman Filter. However, with a trajectory shaping term a different estimator is obtained. The second objective of this work is to investigate the impact and advantages of the trajectory shaping (TS) implementation. The issues under consideration are closely related to the question of the best choice for the disturbance attenuation ratio. To investigate these open issues, two representative models will be used: a Simple Boat Guidance Problem (SBGP) and a Missile Guidance Engagement (MGE) problem. For the former, a closed form analytical solution is obtained. Moreover, to facilitate the understanding of the different factors, a validated approximation to the estimator equation is employed. The second model is more realistic for missile defense, and is treated by numerical experiments, including a comparison with some other guidance methods.
The paper is organized as follows: Section II contains a brief description of the classical LQDG problem with
perfect information case. Section III presents the disturbance attenuation approach. Section IV presents the SBGP and
Section V the MGE problem. Finally, Section VI concludes the paper.

\section{LQDG Problem for Perfect Information}
One of the first formulations of this problem can be found in
\cite{bryson1975applied}. Consider a Linear Time Invariant (LTI) system, where the dynamic equation is given by
\begin{equation} 
\begin{array}{l}
\dot x = Ax + Bu + Dw\\
x\left( {{t_0}} \right) = {x_0}
\end{array}
\end{equation}
where $x \in {R^n}$ is the state vector, $u \in {R^s}$ is the control input signals (the pursuer strategy) and $w \in {R^q}$ is an exogenous disturbance (the evader strategy). The matrices $A,B$ and $D$ are of the appropriate dimensions. The LQDG is a Min-Max problem, which can be defined as follows: Find controls $u $ and $w $ which minimize and maximize (respectively) the cost function
\begin{equation}
\begin{array}{l}
J\left( {u,w} \right) = \frac{1}{2}{x_f}^T{Q_f}{x_f}\\
 + \frac{1}{2}\int\limits_{{t_0}}^{{t_f}} {\left[ {{x^T}Qx + {u^T}u - {\gamma ^2}{w^T}{W^{ - 1}}w} \right]dt} 
\end{array}
\end{equation}
where $Q$ and $W$ are weight matrices and $\gamma^2$ is related to the ratio between the maximizer and minimizer maneuvering capabilities. For convenience, the weight on the control is usually constant and taken as unity. The optimal solution satisfies the "saddle point inequality" \cite{isaacs1999differential}, meaning 
\begin{equation}
J\left( {{u^*},w} \right) \le J^*\left( {{u^*},{w^*}} \right) \le J\left({u,{w^*}} \right)
\end{equation}
where the optimal strategies in the pair $\left\{ {{u^*},{w^*}} \right\}$ are called the saddle point controls. If one player deviates from this strategy, the other player gains advantage (a "zero sum game"). The saddle point value of the cost is defined by

\begin{equation}
{J^*}\left( {{u^*},{w^*}} \right) \buildrel \Delta \over = \frac{1}{2}{x^T}Xx\
\end{equation}

where $X$ will be generated by the following Differential Riccati Equation (DRE)
\begin{equation}
\begin{array}{l}
\dot X + XA + {A^T}X - X\left( {B{B^T} - {\gamma ^{ - 2}}DW{D^T}} \right)X + Q =
0\\
{X^*(t_f)} = {Q_f}.
\end{array}\
\end{equation}
The saddle point inequality is satisfied by the following strategic pair 
\begin{equation}
\begin{array}{l}
{u^*}\left( t \right) =  - {B^T}X(t)x(t)\\
{w^*}\left( t \right) = {\gamma ^{ - 2}}W{D^T}X(t)x(t).
\end{array}
\end{equation}
In the sequel we occasionally omit the dependence on $t$ for the simplicity of notations. In terms of Game Theory, this pair consists of the optimal strategies for both players. Note that when $\gamma  \to \infty $ we can consider the optimal strategy of the evader as zero. Generally, the solution may exist above some critical value, i.e. $\gamma  \ge {\gamma _c}$ which guarantees the existence of the saddle point strategies. Otherwise, if $\gamma  < {\gamma _c}$, there might be no solution for this game. 

\section{LQDG Problem for Imperfect State Information}
The theory is developed in \cite{speyer2010primer}.
Consider noisy continuous measurements   the pursuer. The measurement equation is given by

\begin{equation}
z(t) = Hx(t) + v(t)
\end{equation}

where $H$ is the observation matrix and $v(t)$  is a finite-energy (i.e. square integrable) additive noise (without any statistical prior information). Define a general expression for the system output
\begin{equation}
y = \left[\begin{array}{c} Cx \\ u \end{array}\right]
\end{equation}
where $C$ is of the appropriate dimensions. Dealing with imperfect 
state information involves uncertainty in three factors: evader strategy ($w$), noise strategy  ($v$) and initial condition ($x_0$). Each one of these factors is weighted by positive-definite symmetric matrix: $W, V, Y_0$, respectively. A general representation of the input-output relationship between the disturbances and the system output is defined as the DA ratio
\begin{equation}
{D_a} \buildrel \Delta \over = \frac{\Upsilon }{\Phi }
\end{equation}
where
\begin{equation}
\Upsilon = \frac{1}{2} {{x_f}^T{X_f}{x_f} +
\frac{1}{2}\int\limits_{{t_0}}^{{t_f}} {\left( {{x^T}Qx + {u^T}u} \right)dt} } \
\end{equation}
and
\begin{equation}
\Phi  = \frac{1}{2}x_0^TY_0^{ - 1}{x_0} +\frac{1}{2}\int\limits_{{t_0}}^{{t_f}} {\left( {{w^T}{W^{ - 1}}w + {v^T}{V^{ -1}}v} \right)dt}
\end{equation}
where for simplicity, we define $Q= {C^T}C$; this matrix is the so called "trajectory shaping" matrix.  In order to find the optimal linear strategy of the pursuer $u = u\left( {{Z_t}}\right) $, we use the measurement history which is given by ${Z_t} \equiv \left\{ {z\left( s \right):\,0 \le s \le t} \right\}$ so that the DA ratio is bounded as
\begin{equation}
{D_a} \le \gamma^2 
\end{equation}
for all admissible processes of $w,v$ and $x_0$. Similarly to the perfect
information case, the choice of $\gamma^2$ is not arbitrary. 
There is a critical value ${\gamma_c}^2 $ where if $\gamma^2\le {\gamma_c}^2$, the solution to the problem does not exist.
We consider the following cost function
\begin{equation}
J\left( {u,w} \right) = \Upsilon  - \gamma^2 \Phi.
\end{equation}
After substitution and arranging terms, the cost function is given by 
\begin{equation}
\begin{array}{l}
J = \frac{1}{2}{x_f}^T{X_f}{x_f} - \frac{1}{2}\gamma _{}^2x_0^TY_0^{ - 1}{x_0}\\
 + \frac{1}{2}\int\limits_{{t_0}}^{{t_f}} {\left[ {{x^T}Qx + {u^T}u - \gamma _{}^2\left( {{w^T}{W^{ - 1}}w + {v^T}{V^{ - 1}}v} \right)} \right]dt} .
\end{array}
\end{equation}
The goal in this differential game is to find the Min-Max solution,
\begin{equation}
\begin{array}{l}
{J^*}\left( {{u^*},{{w}^*}} \right) = \mathop {\min }\limits_u \mathop
{\max }\limits_{w} J\left( {u,w} \right)\\
\end{array}
\end{equation}
The optimal DA strategy (control) for the pursuer is becomes:
\begin{equation}
\begin{array}{l}
{u^*} =  - {B^T}X{\left( {I - {\gamma ^{ - 2}}YX} \right)^{ - 1}}\hat x \\
{u^*} =\Lambda\hat x
\end{array}
\end{equation}
where the estimated state $\hat x$ and the measure of the uncertainty 
of the state estimation error $Y$, satisfy the following 
differential equations 
\begin{equation}
\begin{array}{l}
\dot {\hat x} = A\hat x + Bu + {\gamma ^{ - 2}}YQ\hat x + Y{H^T}{V^{ - 1}}\left({z - H\hat x} \right)\\
{{\hat x}_0} = 0
\end{array}
\end{equation}
and,
\begin{equation}
\begin{array}{l}
\dot Y = AY + Y{A^T} + DW{D^T} - Y\left( {{H^T}{V^{ - 1}}H - {\gamma ^{ - 2}}Q} \right)Y\\
Y\left( 0 \right) = {Y_0}.
\end{array}
\end{equation}
One can notice that both Differential Ricatti equations (DRE) as well as
the state estimation equation involve the TS matrix $Q$. 
However, when TS is not considered (i.e. $Q=0$), 
or when $\gamma$ tends to infinity, the familiar (continuous-time) Kalman filter is recovered, 
\begin{equation}
\begin{array}{l}
\dot {\hat x} = A\hat x + Bu + Y{H^T}{V^{ - 1}}\left( {z - H\hat x} \right)\\
{{\hat x}_0} = 0
\end{array}
\end{equation}  
and,
\begin{equation}
\begin{array}{l}
\dot Y = AY + Y{A^T} + DW{D^T} - Y{H^T}{V^{ - 1}}HY\\
Y\left( 0 \right) = {Y_0}.
\end{array}
\end{equation}
The solution  exists if and only if ${\left( {I - {\gamma ^{ - 2}}YX} \right)}$ is positive definite:
\begin{equation}
\Omega  \buildrel \Delta \over = I - {\gamma ^{ - 2}}YX > 0
\end{equation} 
$\forall t \in \left[ {{t_0},{t_f}} \right]$.  Additionally, two positivity conditions should be satisfied, $X > 0,\,\,\,\,Y > 0$.
We next note that the in the solution above of the imperfect state information game, the same guidance gain of (6) of perfect information game, is applied on $\Omega^{-1} \hat x$, where $\hat x$ is just the Kalman filter, when either $Q=0$ (i.e. no trajectory shaping) or $\gamma$ that tends to infinity (i.e. one sided optimal control rather than a game). At this point we note that the above solution to the imperfect state-information game, is just one realization of the pursuer optimal strategy, where different realizations can be readily obtained by applying similarity transformations. Indeed, in \cite{Or2018} one such particularly meaningful transformation has been applied where it has been shown that the similarity transformation utilizing $\Omega^{-1} >0$  leads to the state-estimation $\bar x = \Omega^{-1} \hat x$ which is just the $H_\infty$ filter (see \cite{shaked1990estimation} and \cite{green1994robust}) and where $Y$ of (20) is just the RDE associated with this filter, regardless of the value of $Q$ and also for finite $\gamma$.      	
\section{ A Simple Boat Guidance Problem}
The following problem is based on the modeling used in the famous “Zermelo problem” (Fig. 1). \cite{zermelo1931navigationsproblem}. Two
boats without any dynamics (i.e. of zeroth order) play pursuit- evasion game with imperfect information. We assume full-state measurements to the pursuer, which are, however, corrupted with an additive noise. This so-called simple boat guidance problem (SBGP) was
first considered in \cite{Or2018},\cite{or2018t}.  
\begin{figure}[h]
\includegraphics[width=8.5cm]{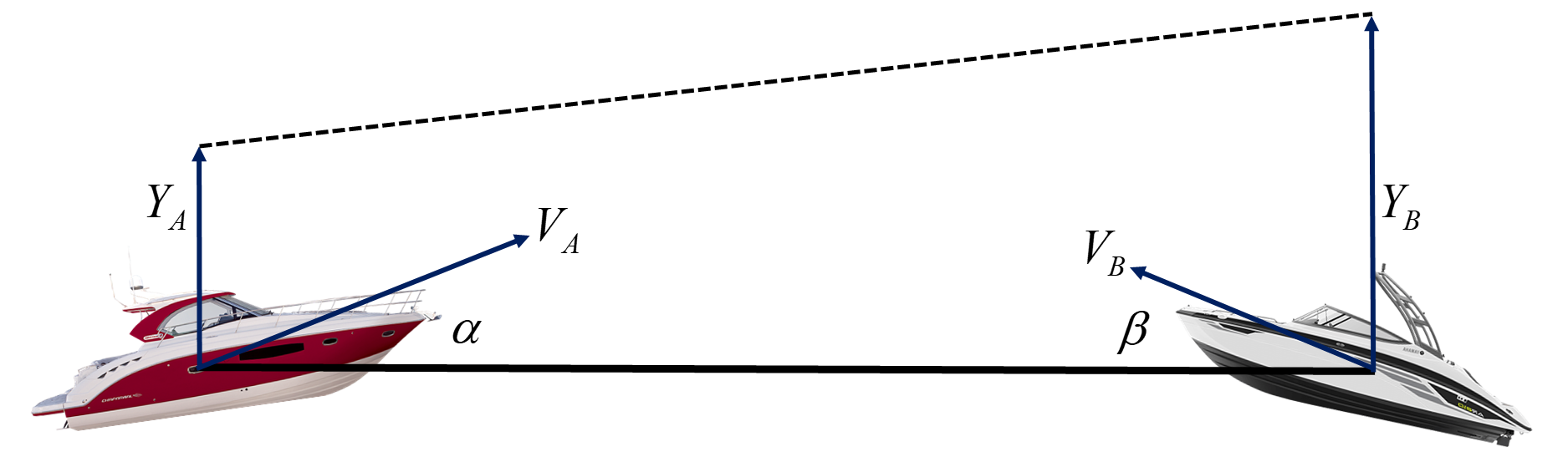}
\centering
\caption{The Simple Boat Guidance Problem (SBGP) Geometry}
\end{figure}
	
\subsection{Scenario Description }
Consider two boats, $A$ and $B$, where $A$ (pursuer) aims to hit boat $B$ (evader). In order to accomplish this mission, boat $A$ controls its heading angle $\alpha$ and tries to navigate towards boat $B$. On the other hand, boat $B$ tries to evade from boat $A$, and does that by controlling tits own heading angle $\beta$. For simplicity,
the following assumptions are taken: 1. A two dimensional problem. 2. Both boats have constant velocities, where boat $A$ is faster than boat $B$ ($V_A>V_B$). 3. The boats have direct control over their heading angles $\alpha$ and $\beta$. 4. Linearized scenario, where $\alpha \ll 1$ , $\beta \ll 1$. The equation of motion is expressed in terms of
\begin{equation}
\begin{array}{l}
x \buildrel \Delta \over = {Y_B} - {Y_A}
\end{array}
\end{equation}
and is described by
\begin{equation}
\dot x = {V_B}\sin \beta  - {V_A}\sin \alpha  \buildrel \Delta \over = w + u
\end{equation}
The measurement equation (for $A$) is
\begin{equation}
z = x + v
\end{equation}
where $x$, the so-called "relative separation" is the difference between the distances ($Y_A, Y_B$) of the boats from the initial Line Of Sight (LOS). $z$ is the measurement performed by boat $A$ and $v \in R^1$ is the additive noise. $u,w \in R^1$ are the controls of the pursuer and the evader, respectively. The problem is formulated as a Min-Max zero-sum game, where $A$ aims to minimize the relative separation at the terminal time, whereas $B$ aims to maximize it.

\subsection{Perfect Information Game}
We first consider the perfect information game where there is no noisy measurement ($v = 0)$. The cost function to be minimzed by $u$ and maximized by $w$, is given by $J$ of (14) with $Q=0$  
\begin{equation}
J\left( {u,w} \right) =\frac{b}{2}{x_f}^2 + \frac{1}{2}\int\limits_0^{{t_f}} {\left( {{u^2} - {\gamma^2}{w^2}} \right)dt} 
\end{equation}
For this simple game, the RDE of (5) reads	
\begin{equation}
\begin{array}{l}
\dot X = \left( {1 - {\gamma ^{ - 2}}} \right){X^2}\,\,\,\,\\
{X(t_f)} = b
\end{array}
\end{equation}
where $b$ is penalization weight on the relative separation at the terminal time. The optimal solution is as follows (see also \cite{ben1998advances})
\begin{equation}
{u^*} =  - \frac{1}{{1/b + \left( {1 - {\gamma ^{ - 2}}} \right)\left( {{t_f} -t} \right)}}x
\end{equation} 
and
\begin{equation}
{w^*} = {\gamma ^{ - 2}}\frac{1}{{1/b + \left( {1 - {\gamma ^{ - 2}}}\right)\left( {{t_f} - t} \right)}}x.
\end{equation}
In the limiting case, where $\gamma  \to \infty$, the evader's actions tend to zero as could be expected from scrutiny of the cost function $J$. Moreover, by letting $\gamma  \to \infty$ and $b  \to \infty$, we obtain a simple Collision Course Guidance (CCG) Law 
\begin{equation}
{u^*} =  - \frac{1}{{\left( {{t_f} - t} \right)}}x
\end{equation}
and
\begin{equation}
{w^*} = 0.
\end{equation}
	
\subsection{Imperfect Information Game}
Consider now an imperfect information game, where additive noise,$v \ne 0$, corrupts the pursuer's measurements:
\begin{equation}
z = x + v.
\end{equation}
The noise term in the cost function is weighted by $V^{-1}$. Additionally, there is uncertainty in the initial condition that is weighted by $Y_0^{-1}$. The modified cost function cost function  becomes:
\begin{equation}
\begin{array}{l}
\mathop {\max }\limits_w \mathop {\min }\limits_u J =  - \frac{1}{2}\gamma _{}^2{Y_0}^{ - 1}{x_0}^2 + \frac{1}{2}b{x_f}^2\\
 + \frac{1}{2}\int\limits_{{t_0}}^{{t_f}} {\left( {{u^2} - \gamma _{}^2\left( {{w^2} + {V^{ - 1}}{v^2}} \right)} \right)dt}. 
\end{array}
\end{equation}    
The optimal pursuer strategy of (16) is then given by
\begin{equation}
{u^*} =  - X{\left( {1 - {\gamma ^{ - 2}}YX} \right)^{ - 1}}\hat x =  -
\frac{1}{{{X^{ - 1}} - {\gamma ^{ - 2}}Y}}\hat x = \Lambda \hat x.
\end{equation}
The estimated state ${\hat{x}}$ is given by a standard Kalman Filter :
\begin{equation}
\begin{array}{l}
\dot {\hat x}= u + Y{V^{ - 1}}\left( {z - \hat x} \right)\\
\hat {x}_0 = 0
\end{array}
\end{equation}
The corresponding RDE of $Y$ is given by
\begin{equation}
\begin{array}{l}
\dot Y = 1 - {Y^2}{V^{ - 1}}\\
Y\left( 0 \right) = {Y_0}
\end{array}
\end{equation}
where the solution can be found in \cite{bar2004estimation}
\begin{equation}
Y = \sqrt V  + \mu ,\,\,\,\mu  \buildrel \Delta \over = \frac{{2\sqrt V }}{{\frac{{{Y_0} + \sqrt V }}{{{Y_0} - \sqrt V }}{e^{\frac{{2t}}{{\sqrt V }}}} - 1}}.
\end{equation}
Note that the solution of RDE for $X$ is identical to the one of the perfect information case. Hence, the optimal DA strategy for the pursuer is given explicitly by
\begin{equation}
{u^*} = \frac{{ - 1}}{{\left( {1 - {\gamma ^{ - 2}}} \right)\left( {{t_f} - t} \right) + 1/b - {\gamma ^{ - 2}}(\sqrt V + \mu) }}\hat x
\end{equation}  
As already stated, there are three conditions in order to guarantee the existence of the above optimal strategies : $X$ must be positive definite, namely, 
\begin{equation}
\left( {1 - {\gamma ^{ - 2}}} \right)\left( {{t_f} - t} \right) + 1/b > 0
\end{equation} 
$Y$  must be positive definite, reading
\begin{equation}
\sqrt V  + \mu > 0
\end{equation}
and finally (to avoid singularity in the control equation),  
\begin{equation}
\left( {1 - {\gamma ^{ - 2}}} \right)\left( {{t_f} - t} \right) + 1/b - {\gamma ^{ - 2}}(\sqrt V + \mu)  >0.
\end{equation}
Notice that the second condition is always satisfied, hence satisfying the third condition entails the satisfaction of the first one.  Reorganizing the last inequality, we get the following inequality involving $\gamma^2$ to be satisfied for all $t \in [t_0,t_f]$ 
\begin{equation}\label{eq:SRC}
\gamma ^2 (1/b + (t_f - t)) > \sqrt V  + \mu + (t_f - t) \end{equation}
\subsection{Approximated Optimal Strategies}
In order to study the behavior of the players, we approximate the optimal strategies. To this end, we consider (following Ref. \cite{bar2004estimation}), the steady-state expression $Y=\sqrt{V}$. The optimal control strategy and estimation equation are simplified to
\begin{equation}
\begin{array}{l}
u = \tilde \Lambda \hat x\\	\tilde \Lambda  = \frac{{ - 1}}{{\left( {1 - {\gamma ^{ - 2}}} \right)\left({{t_f} - t} \right) + 1/b - {\gamma ^{ - 2}}\sqrt V }}
\end{array}
\end{equation}
and
\begin{equation}
\dot {\hat x} = \left( {\tilde \Lambda  - \frac{1}{{\sqrt V }}} \right)\hat x +\frac{1}{{\sqrt V }}z.
\end{equation}
The justification for this approximation is based on the fact that, typically (i.e. for small enough $\sqrt V$) , the exponential term in the control equation diverges quite rapidly. Let the critical DA ratio $\gamma_c^2$ be the minimum value of $\gamma^2$ such that the solution to our problem exists, satisfying the three optimal conditions. The approximate gain $\tilde \Lambda$ must hold $\forall t \in \left[ {{t_0},{t_f}} \right]$, and in particular for $t_f$. By substituting $t=t_f$, the approximated gain approaches to 
\begin{equation}
\tilde \Lambda \left( {{t_f}} \right) =  - \frac{1}{{1/b - {\gamma ^{ - 2}}\sqrt V }}.
\end{equation}
In order to meet the optimality conditions at $t_f$, the following inequality must hold
\begin{equation}
1/b - {\gamma ^{ - 2}}\sqrt V  > 0.
\end{equation}  
Notice that if it does not hold, we have a finite escape time in the obtained optimal strategy. From the last inequality we readily get a lower bound for $\gamma^2$
\begin{equation}
{\gamma ^2} > b\sqrt V. 
\end{equation}
However, we must consider the whole time interval (not only the final time), thus $\gamma^2$ must satisfy the following inequality for $t \in [0,t_f]$ 
\begin{equation}
\begin{array}{l}
\gamma ^2 > \frac{{{t_{go}} + \sqrt V }}{{{t_{go}} + 1/b}}
\buildrel \Delta \over = \psi \left( {{t_{go}}} \right)\\
{t_{go}} \buildrel \Delta \over = {t_f} - t.
\end{array}
\end{equation}
The limiting values of $\gamma_c$ is then obtained combining 
\begin{equation}
\begin{array}{l}
\mathop {\lim }\limits_{{t_{go}} \to t_{f}} \psi \left( {{t_{go}}} \right) = 
\frac{{1 + \sqrt V /{t_{f}}}}{{1 + 1/\left( {b{t_{f}}} \right)}} \\
\mathop {\lim }\limits_{{t_{go}} \to 0} \psi \left( {{t_{go}}} \right) = b\sqrt V. 
\end{array}
\end{equation} 
readily leading to 
\begin{equation}
{\gamma _c}^2 = \max \left\{ {b\sqrt V ,\frac{{1 + \sqrt V /{t_{f}}}}{{1 +1/\left( {b{t_{f}}} \right)}}} \right\}.
\end{equation}
If we want it to hold for all $t_f$, we require
\begin{equation}
{\gamma _c}^2 = \max \left\{ {b\sqrt V ,1} \right\}.
\end{equation}
Under this approximation, we can examine the control behavior with respect to the noise term for a fixed ${\gamma ^2} > \gamma _c^2$. By that we mean that we fix $\gamma$ and solve the problem for different values of the noise term $V$. 
The gain magnitude is a function of $\gamma$, $b$ and $V$: 
\begin{equation}
\left| {\tilde \Lambda } \right| = \frac{1}{{\left( {1 - {\gamma ^{ - 2}}} \right)\left( {{t_f} - t} \right) + 1/b - {\gamma ^{ - 2}}\sqrt V }}.
\end{equation}
One can notice that as the noise term gets higher, the difference between the positive  and negative  values in the denominator reduces. As a result, the gain grows. On the other hand, using the (approximated) critical value of the DA
ratio, adjusting it to the varying $V$ by setting $ \gamma _c^2 = b\sqrt V $ (assuming that this is the dominant term for the lower limit) we get the gain  magnitude a function of $b$ and $V$ as follows: 
\begin{equation}
\left| {\tilde \Lambda } \right| = \frac{1}{{\left( {1 - 1/b\sqrt V } \right)\left( {{t_f} - t} \right)}}.
\end{equation}
Here, the opposite behavior of the gain magnitude with respect to the noise term is evident: as the noise term gets higher, the gain gets lower. 
\section{Missile Guidance Engagement (MGE)}
Consider a guidance engagement problem between two adversaries: a missile and a target. We will employ the following standard assumptions\cite{ben1998advances}:
1.	The engagement scenario is two dimensional, in the horizontal plane. 2. The speeds of the adversaries are assumed to be constant during the engagement. 3.The trajectories can be linearized around a collision triangle. 4. The missile is more maneuverable than the target. 5. The dynamics for both players is assumed to be of first order. Following the geometry in Figure 4, we define: ${x_1} \buildrel \Delta \over = \Delta y$. Let $n_T$  and $a_M$ be the target and missile accelerations
perpendicular to the initial Line-Of-Sight (LOS), and let $w$ and $u$ be the corresponding acceleration commands. The following state-space equations are obtained:

\begin{figure}[h]
\includegraphics[width=8cm]{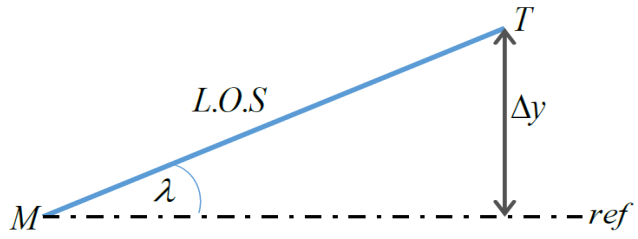}
\centering
\caption{MGE, Problem Geometry}
\end{figure} 
	
\begin{equation}
\begin{array}{l}
\left[ {\begin{array}{*{20}{c}}
{{{\dot x}_1}}\\
{{{\dot x}_2}}\\
{{{\dot n}_T}}\\
{{{\dot a}_M}}
\end{array}} \right] = \left[ {\begin{array}{*{20}{c}}
0&1&0&0\\
0&0&1&{ - 1}\\
0&0&{ - {\theta ^{ - 1}}}&0\\
0&0&0&{ - {T^{ - 1}}}
\end{array}} \right]\left[ {\begin{array}{*{20}{c}}
{{x_1}}\\
{{x_2}}\\
{{n_T}}\\
{{a_M}}
\end{array}} \right]\\
\,\,\,\,\,\,\,\,\,\,\,\,\, + \left[ {\begin{array}{*{20}{c}}
0\\
0\\
0\\
{{T^{ - 1}}}
\end{array}} \right]u + \left[ {\begin{array}{*{20}{c}}
0\\
0\\
{{\theta ^{ - 1}}}\\
0
\end{array}} \right]w
\end{array}
	\end{equation}
	
	where $T$ and $\theta$ are the time constants of the target and the missile
	respectively. 
	Define
	
	\begin{equation}
	\begin{array}{l}
	{{x}} \buildrel \Delta \over = {\left[ {\begin{array}{*{20}{c}}
			{{x_1}}&{{x_2}}&{{n_T}}&{{a_M}}
			\end{array}} \right]^T}\\
	{{B}} \buildrel \Delta \over = {\left[ {\begin{array}{*{20}{c}}
			0&0&0&{{T^{ - 1}}}
			\end{array}} \right]^T}\\
	{{D}} \buildrel \Delta \over = {\left[ {\begin{array}{*{20}{c}}
			0&0&{{\theta ^{ - 1}}}&0
			\end{array}} \right]^T}\\
	{{A}} \buildrel \Delta \over = \left[ {\begin{array}{*{20}{c}}
		0&1&0&0\\
		0&0&1&{ - 1}\\
		0&0&{ - {\theta ^{ - 1}}}&0\\
		0&0&0&{ - {T^{ - 1}}}
		\end{array}} \right]
	\end{array}
	\end{equation}
	Thus 
	\begin{equation}
	\begin{array}{l}
	{{\dot x = Ax + Bu + Dw}}\\
	{{x}}\left( {{0}} \right){\bf{ = }}{{{x}}_{{0}}}
	\end{array}
	\end{equation}
Let the measurement be the LOS angle, $\lambda =  {\Delta y / R}$, where $R$ is
	the range between the adversaries. Hence 
\begin{equation}
\begin{array}{l}
z = Hx + v\\
H = \left[ {\begin{array}{*{20}{c}}
{{{\left( {{V_c}{t_{go}}} \right)}^{ - 1}}}&0&0&0
\end{array}} \right]
\end{array}
\end{equation}
$V_c$ is the closing velocity and $t_{go}=t_f-t$. The cost function is, as before:
\begin{equation}
\begin{array}{l}
J = \frac{1}{2}{x_f}^T{X_f}{x_f} - \frac{1}{2}\gamma _{}^2x_0^TY_0^{ - 1}{x_0}\\
 + \frac{1}{2}\int\limits_{{t_0}}^{{t_f}} {\left[ {{u^T}u - \gamma _{}^2\left( {{w^T}{W^{ - 1}}w + {v^T}{V^{ - 1}}v} \right)} \right]dt} 
\end{array}
\end{equation}
where $u$ is a minimizer and $w,v,x_0$ are maximizers of $J$. \\
As opposed to the previous simpler problem, analytical solutions were not sought after, and a numerical study was performed. 
To this end the following parameters have been used: ${X_f} = diag\left\{
{b,0,0,0} \right\}, {V_c} = 300\left[ {m/s} \right], {x_0} = 0,{Y_0} = {I_4},W = 3,b = 1000, T = 0.1\left[ s \right],\theta  = 0.5\left[ s \right],\,\eta  =\sqrt{V}= 0.5 \times 10^{-3}$, $ w=1g $. The value of $V$ will vary according the applied noise variance (\cite{ben2017games}). 
In this numerical study we have three objectives:   
1. Studying the behavior of the control gain with fixed and critical DA ratios.  \\
2. Comparing between the DA control, a separation-based control, the perfect	information control and proportional navigation (PN). \\
3. Studying the effects of trajectory shaping on the solution. \\

\subsubsection{DA ratio (Fixed and Critical)}
One advantage of using the LQDG formulation for missile guidance is that one gets Proportional Navigation (PN) like solutions with a continuous control function. The general expression for PN strategy is given by 
\begin{equation}
u  = N'{V_c}\dot \lambda.
\end{equation}
Taking the time derivative of $\lambda$ leads to  
\begin{equation}
u = N'{V_c}\left[ {{x_1}\frac{1}{{{V_c}t_{go}^2}} +{x_2}\frac{1}{{{V_c}t_{go}^{}}}} \right] = \frac{1}{{t_{go}^2}}N'{x_1} +\frac{1}{{t_{go}^{}}}N'{x_2}.
\end{equation}
Hence one can retrieve the equivalent gain $N'$ from the  differential game solution, by identifying
\begin{equation}
N' = {\Lambda _{{x_1}}}t_{go}^2  
\end{equation}
where $ {\Lambda _{{x_1}}} $ is the feedback gain from the state $x_1$ to the optimal control. We aim here at understanding the behavior of this gain under fixed and (near) critical DA strategies. To this end, numerical simulations with different noise levels,$v \sim {\cal N}\left( {0,{{10}^{ - 3}}\eta } \right)$ , were used.
Using the true critical DA ratio (${\gamma _c}$ ) is not advisable in practical applications or numerical simulations.  
Hence the distance from singularity was calculated by the  determinant $\left| \Omega  \right| $ of $\Omega = {\left( {I - {\gamma ^{ - 2}}YX} \right)}$, requiring it to be above a certain threshold e.g. $\left| \Omega  \right| > 0.36$. The optimal control gains for (near) critical DA ratio and fixed DA ratio are depicted in Fig.3 and Fig.4, respectively.The corresponding $\left| \Omega  \right|$'s are depicted in Fig. 5 and Fig. 6. The main conclusion is that by considering the critical DA ratio (which depends on the noise level), we get that the equivalent navigation constant decreases when noise level grows. On the other hand, when we consider some fixed DA ratio (a fixed value for all noise levels) the equivalent navigation constant  grows with the noise level. This behavior is in complete agreement with the analytic results of the previous section. 
Table 1 compares the performance of the fixed versus (near) critical $ \gamma $. $\gamma_c=3 $ for the case $ \eta=0.9 $ (middle row). For a lower values of noise $ \eta=0.6 $ and $ \eta=0.5 $, the corresponding $ \gamma_c=2.55 $ and $ \gamma_c=2.39 $ have been obtained. Evidently, using these (near) critical values improve the performance, both in  miss distances and control efforts.       
 \begin{figure}[!h]
\includegraphics[width=8cm]{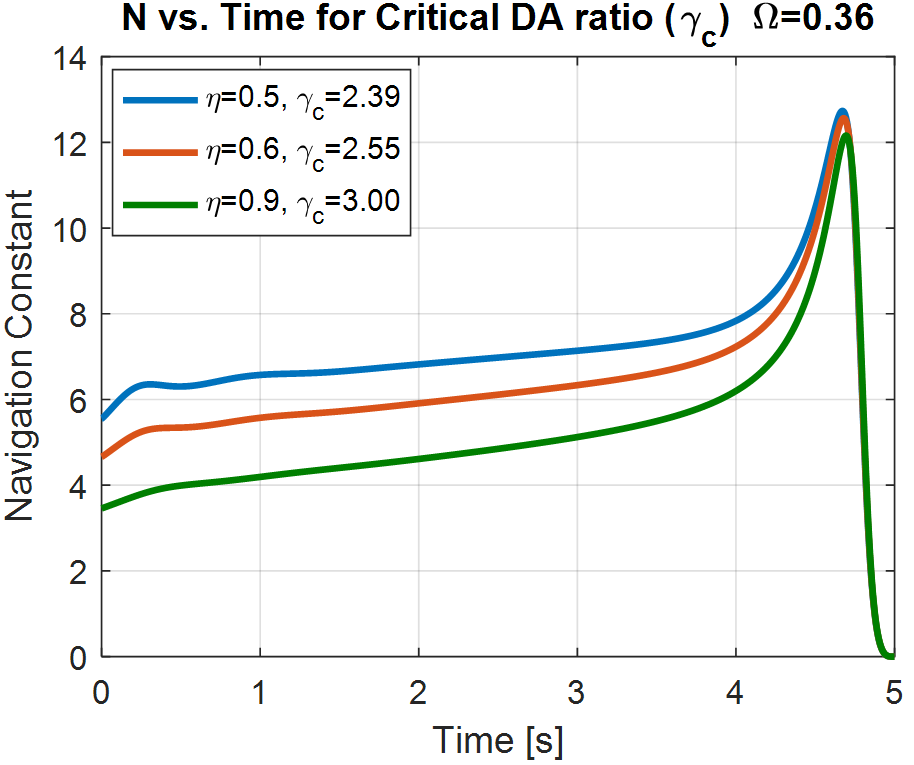}
\centering
\caption{N Constant for Critical DA Ratio}
\end{figure} 

 \begin{figure}[!h]
\includegraphics[width=8cm]{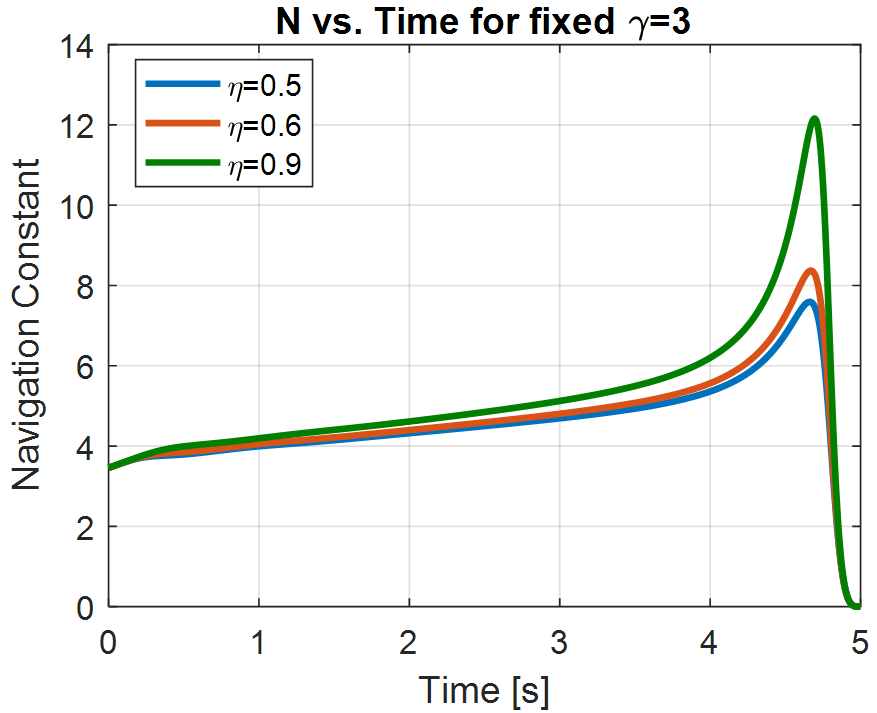}
\centering
\caption{N Constant for Fixed DA Ratio}
\end{figure} 	

 \begin{figure}[!h]
\includegraphics[width=8cm]{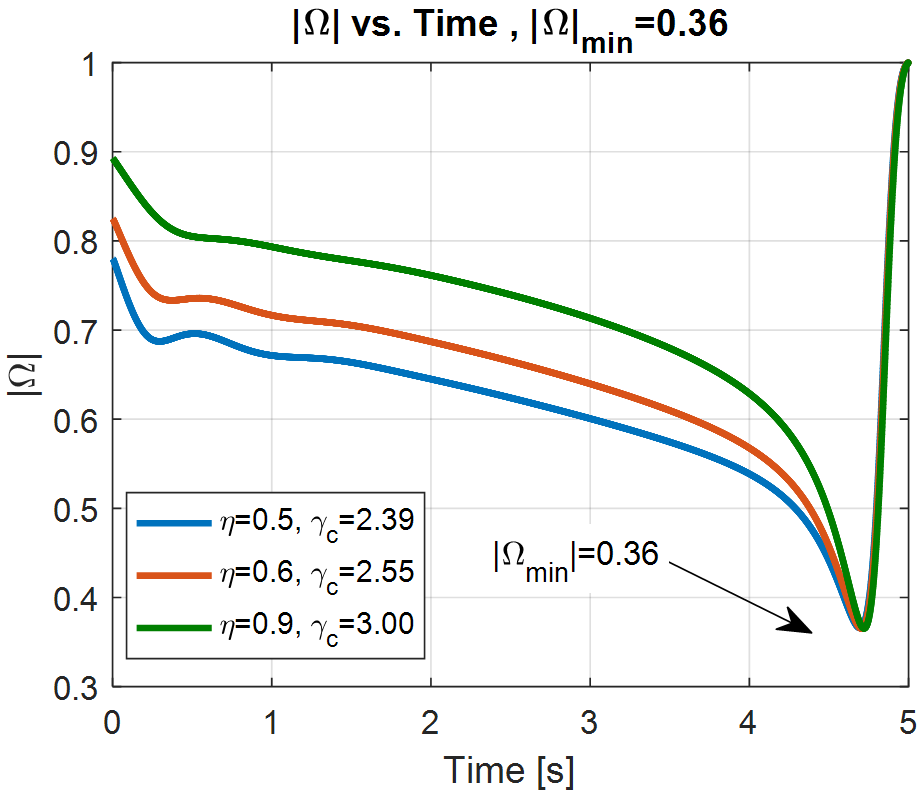}
\centering
\caption{Omega vs. Time (Critical DA)}
\end{figure} 

 \begin{figure}[!h]
\includegraphics[width=8cm]{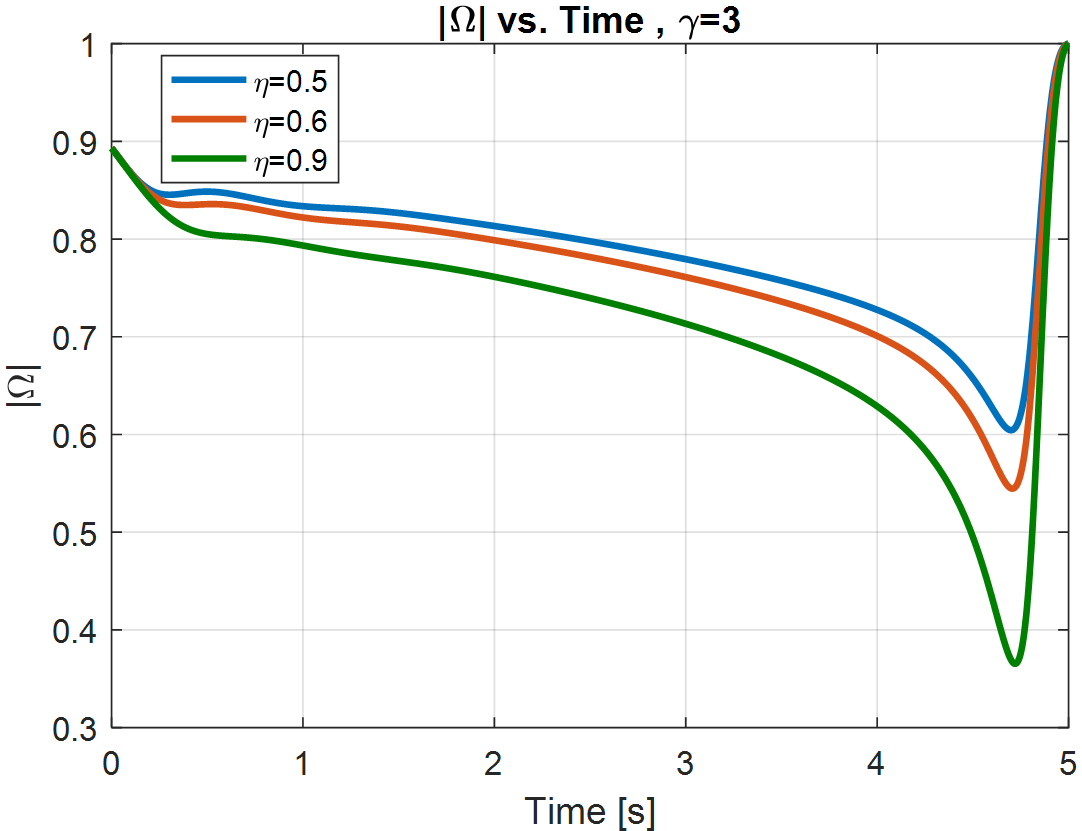}
\centering
\caption{Omega vs. Time (Fixed DA)}
\end{figure}

	\begin{table}[hbt!]
	\caption{\label{tab:table1} MGE Simulation results}
	\centering
	\begin{tabular}{lcccccc}
	\hline
Control  & CEP[cm] & Control Effort \\\hline
	$\gamma  = 3.00, \eta  = 0.5$  &  25    &    556\\
	$\gamma  = 3.00, \eta  = 0.6$  &  26    &    568\\
	$\gamma  = 3.00, \eta  = 0.9$  &  27    &    602\\
	$\gamma_c= 2.39, \eta  = 0.5$  &  20    &    436\\
	$\gamma_c= 2.55, \eta  = 0.6$  &  20    &    481\\
	\hline
	\end{tabular}
	\end{table}
\subsubsection{Comparison between different control strategies}
Consider four controls strategies as given in Table II. In our simulations $u_i$ is restricted to ${\left| {{u_i}} \right| < u_{sat}^{}=4g}$  to avoid  saturation. $\gamma  = 2.5$ was employed. Table III summarizes the performance results. Evidently, DA guidance outperforms PN and the separation-based control, and is close to the perfect information solution.  
\begin{table}[!h]
\caption{\label{tab:table1} MGE Controls}
\centering
\begin{tabular}{lcccccc}
\hline
Name & Equation \\\hline
DA control  &  ${u_1} =  - {B^T}X{\Omega ^{ - 1}}\hat x$ \\
Perfect State Control   & ${u_2} =  - {B^T}Xx$ \\
Separation Control & ${u_3} =  - {B^T}X\hat x$ \\
Proportional Navigation Control & ${u_4} =  - \frac{3}{{t_{go}^2}}\left( {{{\hat x}_1} + t_{go}^{}{{\hat x}_2}} \right)$ \\
\hline
\end{tabular}
\end{table}
\begin{table}[!h]
\caption{\label{tab:table2} MGE Simulation results}
\centering
\begin{tabular}{lcccccc}
\hline
Control  & CEP[cm] & Control Effort \\\hline
$u_1$  &  30    &    392\\
$u_2$  &  14    &    284\\
$u_3$  &  52    &    532\\
$u_4$  &  128    &    960\\
\hline
\end{tabular}
\end{table}
	
\subsubsection{Trajectory shaping}
An advantage of the LQDG formulation is its flexibility, which enables it not only to include in the cost function additional weights on other terminal variables, but also to introduce some Trajectory Shaping by augmenting the cost function with a running cost term on the state variable. The cost function with Q term becomes 
\begin{equation}
\begin{array}{l}
J = \frac{1}{2}{x_f}^T{Q_f}{x_f} - \frac{1}{2}\gamma _{}^2x_0^TY_0^{ - 1}{x_0}\\
 + \frac{1}{2}\int\limits_{{t_0}}^{{t_f}} {\left[ {{x^T}Qx + {u^T}u - \gamma _{}^2\left( {{w^T}{W^{ - 1}}w + {v^T}{V^{ - 1}}v} \right)} \right]dt} 
\end{array}
\end{equation}
where $Q = diag\left\{ {{q_{11}},0,0,0} \right\}$. The following scenario was simulated 	using the parameters: ${V_c} = 300\left[ {m/s}	\right],\,{X_f} = diag\left\{ {b,0,0,0} \right\},\,{x_{i0}} = 0\,\,\,i \in \left\{ {1,3,4} \right\},\,\,{x_{20}} = 5\left[ {m/s} \right],w = 1g,\,\,\,{Y_0}
= {I_4}\,,\,W = 3,\,b = 1000,\,v \sim {\cal N}\left( {0,{{10}^{ - 3}}\eta } \right)$.
The determinant of $ \Omega $ is depicted in Fig. 7 and the corresponding trajectories in Fig. 8. An important observation is that a new local minimum for $|\Omega|$ appears near $t_0$ due to the trajectory shaping term  in addition to the one near $t_f$. Thus, increasing  the trajectory shaping term (without readjusting $\gamma$), may violate the necessary	condition for optimality. It is interesting to note that, by (nearly) equating the two minima, the solutions become similar.
	
\begin{figure}[h]
\includegraphics[width=8cm]{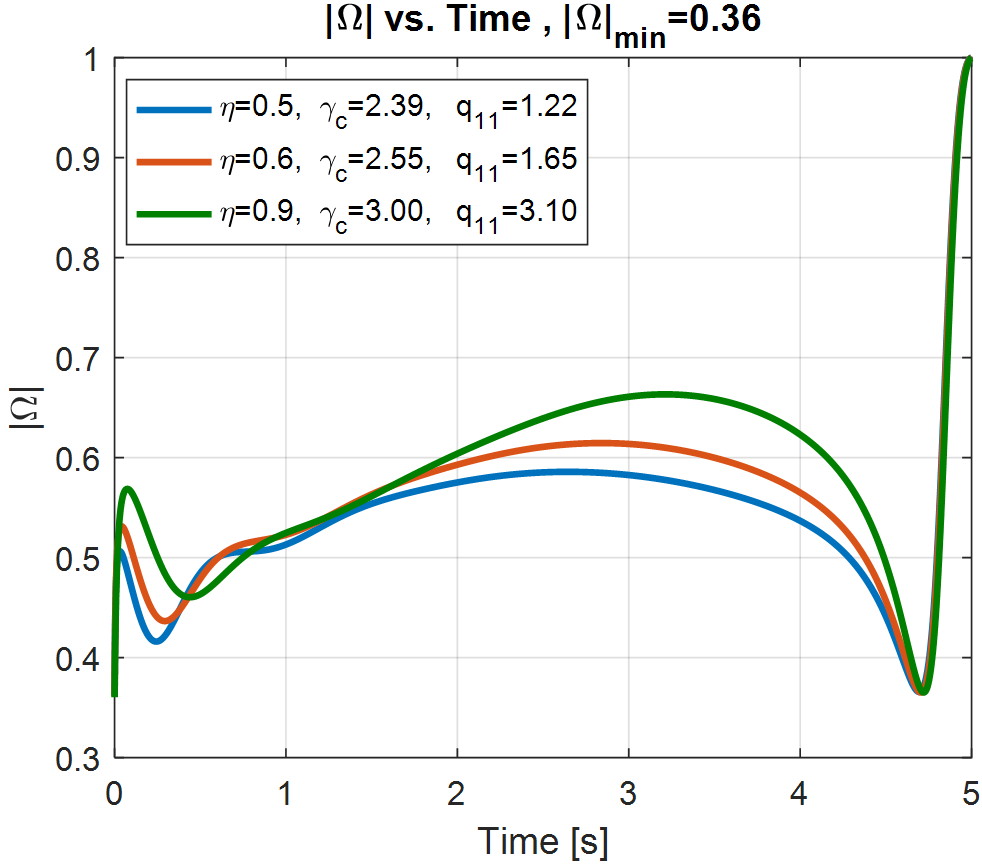}
\centering
\caption{TS - Omega vs. Time for Critical DA ratio}
\end{figure}	

\begin{figure}[h]
\includegraphics[width=8cm]{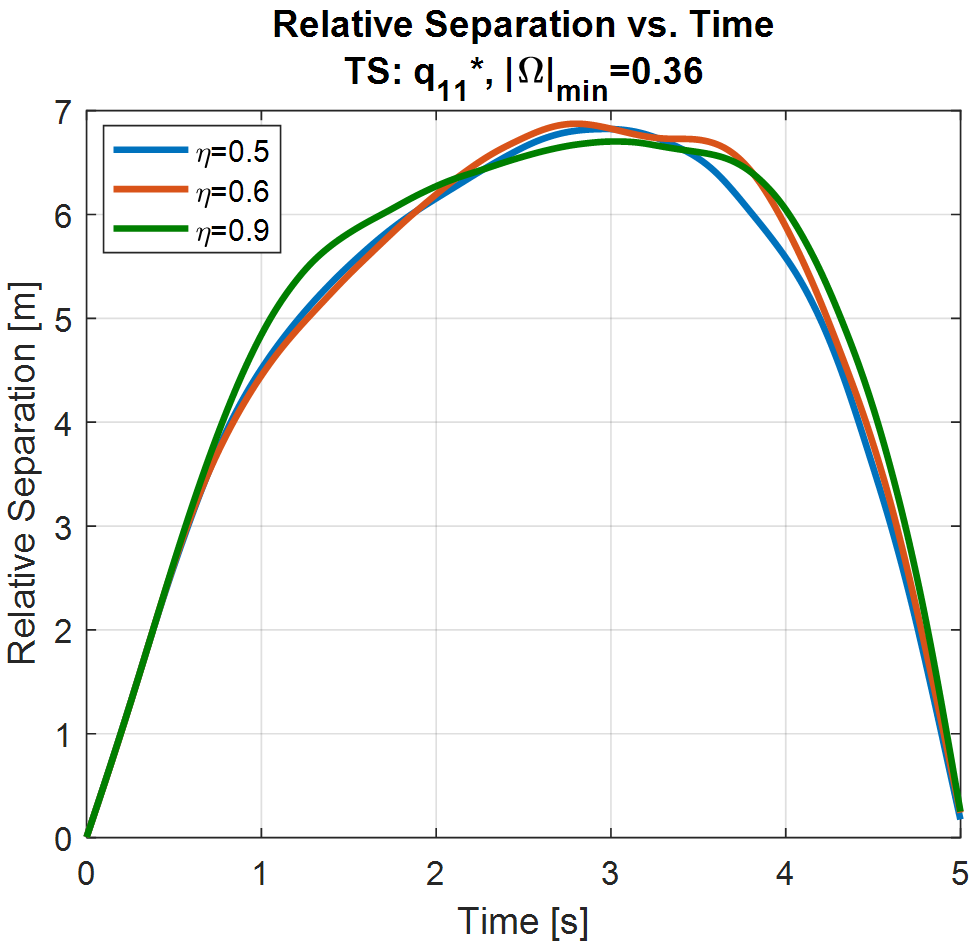}
\centering
\caption{TS - Relative Separation  vs. Time for Critical DA ratio}
\end{figure}
\section{Conclusions}
The problem of DA with imperfect information pattern has been revisited addressing some unanswered questions. Specifically, the issue of the effect of noise on the control strategy, and the effect of trajectory shaping have been studied. A closed form solution was obtained for a representative case – a Simple Boat Guidance Problem. Numerical solutions were obtained for the common guidance problem of Missiles Guidance Engagement. The advantages of using the critical value of the DA ratio over a fixed DA ratio have been demonstrated. Trajectory Shaping was introduced in the MGE game and yielded interesting results regarding two local minima for the spectral radius related term $\Omega$. A comparison of the DA Control, Separated Control and the PN for the MGE was performed. DA Control outperformed both Separated Control and PN. Although this research focused on rather simple guidance problems, it is believed to capture some of the main characteristics of the more general DA guidance approach.
\bibliographystyle{IEEEtran}
\bibliography{Main}
\end{document}